\documentclass[12pt]{amsart}
\usepackage{amssymb}
\usepackage[all]{xy}
\usepackage{tikz}

\usepackage{amscd}
\usepackage{amsmath}
\input xy
\xyoption{all}

\newtheorem{lemma}{Lemma}[section]

\newtheorem{theorem}[lemma]{Theorem}
\newtheorem{proposition}[lemma]{Proposition}

\theoremstyle{definition}
\newtheorem{remark}[lemma]{Remark}
\newtheorem{definition}[lemma]{Definition}

\DeclareMathOperator{\Mod}{Mod}
\DeclareMathOperator{\modd}{mod}

\DeclareMathOperator{\Hom}{Hom}

\DeclareMathOperator{\Ext}{Ext}

\DeclareMathOperator{\Ker}{Ker}
\DeclareMathOperator{\Coker}{Coker}
\DeclareMathOperator{\Imm}{Im}

\newtheorem*{theorem a*}{Theorem A}
\newtheorem*{theorem b*}{Theorem B}

\newcounter{diagram}
\numberwithin{diagram}{section}

\newenvironment{diagram}
  {\stepcounter{diagram}\par\smallskip\noindent\begin{minipage}{\linewidth}\centering}
  {\par Diagram~\thediagram\end{minipage}\par\smallskip}

\begin{document}

\title{higher auslander's formula}

\author{Ramin Ebrahimi}
\address{Department of Pure Mathematics\\
Faculty of Mathematics and Statistics\\
University of Isfahan\\
P.O. Box: 81746-73441, Isfahan, Iran}
\email{ramin69@sci.ui.ac.ir}

\author{Alireza Nasr-Isfahani}
\address{Department of Pure Mathematics\\
Faculty of Mathematics and Statistics\\
University of Isfahan\\
P.O. Box: 81746-73441, Isfahan, Iran\\ and School of Mathematics, Institute for Research in Fundamental Sciences (IPM), P.O. Box: 19395-5746, Tehran, Iran}
\email{nasr$_{-}$a@sci.ui.ac.ir / nasr@ipm.ir}



\begin{abstract}
Let $\mathcal{M}$ be a small $n$-abelian category. We show that the category of finitely presented functors $\modd$-$\mathcal{M}$ modulo the subcategory of effaceable functors $\modd_0$-$\mathcal{M}$ has an $n$-cluster tilting subcategory which is equivalent to $\mathcal{M}$. This  gives a higher-dimensional version of Auslander's formula.
\end{abstract}

\maketitle


\section{Introduction}
For a positive integer $n$, $n$-cluster tilting subcategories of abelian categories were introduced by Iyama \cite{I2, I1} (see also \cite{I3} and \cite{I4}) to construct a higher Auslander correspondence and develop the higher-dimensional analogue of Auslander-Reiten theory. Recently, Jasso \cite{J} introduced $n$-abelian categories as higher-dimensional analogue of abelian categories from the viewpoint of higher homological
algebra. $n$-abelian categories are an axiomatisation of $n$-cluster tilting subcategories. In these subcategories kernels and cokernels don't necessarily exist and are replaced with $n$-kernels and $n$-cokernels. Jasso in \cite{J} proved that any $n$-cluster tilting subcategory of an abelian category is $n$-abelian.

Let $\mathcal{M}$ be a projectively generated small $n$-abelian category, $\mathcal{P}$ be the
subcategory of projective objects in $\mathcal{M}$ and $F:\mathcal{M}\rightarrow \modd$-$\mathcal{P}$ be the functor defined
by $F(X)=\mathcal{M}(-, X)|_{\mathcal{P}}$. Jasso in \cite{J} proved that if there exists an exact duality $D:
\modd$-$\mathcal{P}\rightarrow \modd$-$\mathcal{P}^{op}$, then $F$ is fully faithful and the essential image of $F$ is an $n$-cluster tilting subcategory of $\modd$-$\mathcal{P}$. Kvamme in \cite{Kv} proved that the existence of exact
duality is unnecessary. In this paper we show that the projectively generated assumption is also unnecessary.

Let $\mathcal{M}$ be a small $n$-abelian category. In \cite{EN} we showed that there exist an abelian category $\mathcal{A}$ and a fully faithful functor $\mathcal{M}\hookrightarrow \mathcal{A}$ which preserves and reflects $n$-exact sequences (see Definition \ref{d1}). Motivated by this result it is natural to ask the following question. Does there exist an abelian category $\mathcal{B}$ and a fully faithful functor $\mathcal{M}\hookrightarrow \mathcal{B}$ that $\mathcal{M}$ is equivalent to an $n$-cluster tilting subcategory of an abelian category $\mathcal{B}$ (see \cite[Question 4.6]{EN})? In this paper we give a positive answer to this question. This question is answered independently by Kvamme in \cite{Kv1}.

Let $\mathcal{C}$ be a small abelian category and $\mathcal{G}$ be the category of abelian groups. All contravariant additive functors $F:\mathcal{C}\rightarrow \mathcal{G}$, and natural transformations between them, form an abelian category which is denoted by $\Mod$-$\mathcal{C}$ or $(\mathcal{C}^{op},\mathcal{G})$. A contravariant additive functor $F$ is called {\em finitely presented} if there exist $C_1, C_2\in \mathcal{C}$ and an exact sequence $\Hom_\mathcal{C}(-,C_1)\rightarrow \Hom_\mathcal{C}(-,C_2)\rightarrow F\rightarrow 0$ in $\Mod$-$\mathcal{C}$. The full subcategory of $\Mod$-$\mathcal{C}$ consisting of all finitely presented functors denoted by $\modd$-$\mathcal{C}$. A finitely presented functor $F$ is called {\em effaceable} if there exists a presentation $\Hom_\mathcal{C}(-,C_1)\overset{(-,f)}{\longrightarrow} \Hom_\mathcal{C}(-,C_2)\rightarrow F\rightarrow 0$ such that $f:C_1\rightarrow C_2$ is an epimorphism. The full subcategory of $\modd$-$\mathcal{C}$ consisting of all effaceable functors, denoted by $\modd_0$-$\mathcal{C}$, is a Serre subcategory of $\modd$-$\mathcal{C}$ \cite[Page 202]{Au}. A very famous and important formula due to Auslander \cite[Page 205]{Au}, which is known as Auslander's formula, shows that we have the following equivalence of abelian categories

\begin{equation}
\frac{\rm{mod}\text{-}\mathcal{C}}{\rm{mod}_0\text{-}\mathcal{C}}\simeq \mathcal{C}. \notag
\end{equation}

This equivalence shows that to study the abelian category $\mathcal{C}$ it is enough to study the category of finitely presented functors that has better homological behaviour. This leads to the functorial approach in representation theory, and as Lenzing said in \cite{Le}, a substantial part of Maurice Auslander's work on the representation theory of Artin algebras can be linked to the formula. Krause in \cite{K} proved a derived version of Auslander's formula.

In this paper we prove a higher-dimensional version of Auslander's formula for small $n$-abelian categories. More precisely we prove the following theorem.

\begin{theorem a*}$($Theorem \ref{m2}$)$
Let $\mathcal{M}$ be a small $n$-abelian category. The category of finitely presented functors $\modd$-$\mathcal{M}$ modulo the subcategory of effaceable functors $\modd_0$-$\mathcal{M}$ has an $n$-cluster tilting subcategory which is equivalent to $\mathcal{M}$.
\end{theorem a*}

The paper is organised as follows. In section 2 we recall the definitions of $n$-abelian categories, $n$-cluster tilting subcategories and recall some results that we need in the paper. In Section 3, we first recall the definitions of Serre subcategories and localisation in abelian categories. Then we show that for a small $n$-abelian category $\mathcal{M}$, $\modd_0$-$\mathcal{M}$ is a Serre subcategory of $\modd$-$\mathcal{M}$ and hence we have the full and faithful functor
 \begin{equation}
\widetilde{H}:\mathcal{M}\longrightarrow \frac{\rm{mod}\text{-}\mathcal{M}}{\rm{mod}_0\text{-}\mathcal{M}} \notag
\end{equation}
where $\widetilde{H}$ is a composition of the Yoneda functor $H:\mathcal{M}\rightarrow \rm{mod}\text{-}\mathcal{M}$ by the localisation functor. In Section 4, following Auslander, we describe the category of finitely presented functors on $\mathcal{M}$ modulo the subcategory of effaceable functors and show that $\mathcal{M}$ is equivalent to an $n$-cluster tilting subcategory of this quotient category.

\subsection{Notation}
Throughout this paper, unless otherwise stated, $n$ always denotes a fixed positive integer, $\mathcal{M}$ is a fixed small $n$-abelian category and $\mathcal{G}$ is the abelian category of all abelian groups.


\section{preliminaries}
In this section we recall the definitions of $n$-abelian category, $n$-cluster tilting subcategory and recall some results that we need in the rest of paper. For further information and motivation of definitions the readers are referred to \cite{I2, I3, I1, J}.

\subsection{$n$-abelian categories}
Let $\mathcal{A}$ be an additive category and $f:A\rightarrow B$ a morphism in $\mathcal{A}$. A {\em weak cokernel} of $f$ is a morphism $g:B\rightarrow C$ such that for all $C^{\prime} \in \mathcal{A}$  the sequence of abelian groups
\begin{center}
$\begin{CD}
\Hom_\mathcal{A}(C,C')@ > \Hom_\mathcal{A}(g,C')>> \Hom_\mathcal{A}(B,C')@ > \Hom_\mathcal{A}(f,C')>>  \Hom_\mathcal{A}(A,C'),
\end{CD}$
\end{center}
is exact. The concept of {\em weak kernel} is defined dually.

Let $d^0:X^0 \rightarrow X^1$ be a morphism in $\mathcal{A}$. An {\em $n$-cokernel} of $d^0$ is a sequence
\begin{equation}
(d^1, \ldots, d^n): X^1 \overset{d^1}{\rightarrow} X^2 \overset{d^2}{\rightarrow}\cdots \overset{d^{n-1}}{\rightarrow} X^n \overset{d^n}{\rightarrow} X^{n+1} \notag
\end{equation}
of objects and morphisms in $\mathcal{A}$, such that for each $Y\in \mathcal{A}$
the induced sequence of abelian groups
\begin{align}
0 \rightarrow \Hom_\mathcal{A}(X^{n+1},Y) \rightarrow \Hom_\mathcal{A}(X^n,Y) \rightarrow\cdots\rightarrow \Hom_\mathcal{A}(X^1,Y) \rightarrow \Hom_\mathcal{A}(X^0,Y) \notag
\end{align}
is exact. Equivalently, the sequence $(d^1, \ldots, d^n)$ is an $n$-cokernel of $d^0$ if for all $1\leq k\leq n-1$
the morphism $d^k$ is a weak cokernel of $d^{k-1}$, and $d^n$ is moreover a cokernel of $d^{n-1}$ \cite[Definition 2.2]{J}. The concept of {\em $n$-kernel} of a morphism is defined dually.
\begin{definition}(\cite[Definition 2.4]{L})\label{d1}
Let $\mathcal{A}$ be an additive category. A {\em right $n$-exact sequence} in $\mathcal{A}$ is a complex
\begin{equation}
X^0 \overset{d^0}{\rightarrow} X^1 \overset{d^1}{\rightarrow} \cdots \overset{d^{n-1}}{\rightarrow} X^n \overset{d^n}{\rightarrow} X^{n+1} \notag
\end{equation}
such that $(d^1, \ldots, d^{n})$ is an $n$-cokernel of $d^0$. The concept of {\em left $n$-exact sequence} is defined dually. An {\em $n$-exact sequence} is a sequence which is both a right $n$-exact sequence and a left $n$-exact sequence.
\end{definition}

Let $\mathcal{A}$ be a category and $A$ be an object of $\mathcal{A}$. A morphism $e\in \Hom_\mathcal{A}(A, A)$ is {\em idempotent} if $e^2 = e$. $\mathcal{A}$ is called {\em idempotent complete}
if for every idempotent $e\in \Hom_\mathcal{A}(A, A)$ there exist an object $B$ and morphisms $f\in\Hom_\mathcal{A}(A, B)$ and
$g\in\Hom_\mathcal{A}(B, A)$ such that $gf = e$ and $fg = 1_B$ \cite[Page 61]{Fr}.

\begin{definition}$($\cite[Definition 3.1]{J}$)$
An {\em $n$-abelian} category is an additive category $\mathcal{M}$ which satisfies the following axioms.
\begin{itemize}
\item[(i)]
The category $\mathcal{M}$ is idempotent complete.
\item[(ii)]
Every morphism in $\mathcal{M}$ has an $n$-kernel and an $n$-cokernel.
\item[(iii)]
For every monomorphism $d^0:X^0 \rightarrow X^1$ in $\mathcal{M}$ and for every $n$-cokernel $(d^1, \ldots, d^n)$ of $d^0$, the  following sequence is $n$-exact:
\begin{equation}
X^0 \overset{d^0}{\rightarrow} X^1 \overset{d^1}{\rightarrow} \cdots \overset{d^{n-1}}{\rightarrow} X^n \overset{d^n}{\rightarrow} X^{n+1}. \notag
\end{equation}
\item[(iv)]
For every epimorphism $d^n:X^n \rightarrow X^{n+1}$ in $\mathcal{M}$ and for every $n$-kernel $(d^0, \ldots, d^{n-1})$ of $d^n$, the following sequence is $n$-exact:
\begin{equation}
X^0 \overset{d^0}{\rightarrow} X^1 \overset{d^1}{\rightarrow} \cdots \overset{d^{n-1}}{\rightarrow} X^n \overset{d^n}{\rightarrow} X^{n+1}. \notag
\end{equation}
\end{itemize}
\end{definition}

A subcategory $\mathcal{B}$ of an abelian category $\mathcal{A}$ is called {\em cogenerating} if for every object
$X\in \mathcal{A}$ there exist an object $Y\in\mathcal{B}$ and a monomorphism $X\rightarrow Y$. The concept of
{\em generating} subcategory is defined dually \cite[Page 724]{J}.

Let $\mathcal{A}$ be an additive category and $\mathcal{B}$ be a full subcategory of $\mathcal{A}$. $\mathcal{B}$ is called {\em covariantly finite in $\mathcal{A}$} if for every $A\in \mathcal{A}$ there exist an object $B\in\mathcal{B}$ and a morphism
$f : A\rightarrow B$ such that, for all $B'\in\mathcal{B}$, the sequence of abelian groups $\Hom_\mathcal{A}(B, B')\rightarrow \Hom_\mathcal{A}(A, B')\rightarrow 0$ is exact. Such a morphism $f$ is called a {\em left $\mathcal{B}$-approximation of $A$}. The notions of {\em contravariantly
finite subcategory of $\mathcal{A}$} and {\em right $\mathcal{B}$-approximation} are defined dually. A {\em functorially
finite subcategory of $\mathcal{A}$} is a subcategory which is both covariantly and contravariantly finite
in $\mathcal{A}$ \cite[Page 113]{AR}.

\begin{definition}$($\cite[Definition 3.14]{J}$)$
Let $\mathcal{A}$ be an abelian category and $\mathcal{M}$ be a generating-cogenerating full subcategory of $\mathcal{A}$. $\mathcal{M}$ is called an {\em $n$-cluster tilting subcategory of $\mathcal{A}$} if $\mathcal{M}$ is functorially finite in $\mathcal{A}$ and
\begin{align}
\mathcal{M}& = \{ X\in \mathcal{A} \mid \forall i\in \{1, \ldots, n-1 \}, \Ext_{\mathcal{A}}^i(X,\mathcal{M})=0 \}\notag \\
& =\{ X\in \mathcal{A} \mid \forall i\in \{1, \ldots, n-1 \}, \Ext_{\mathcal{A}}^i(\mathcal{M},X)=0 \}.\notag
\end{align}

Note that $\mathcal{A}$ itself is the unique $1$-cluster tilting subcategory of $\mathcal{A}$.

\end{definition}
\begin{remark}$($\cite[Remark 3.15]{J}$)$ {\rm Let $\mathcal{A}$ be an abelian category and $\mathcal{M}$ be an $n$-cluster tilting subcategory of $\mathcal{A}$. Since $\mathcal{M}$ is a generating-cogenerating subcategory of $\mathcal{A}$, for each $A\in \mathcal{A}$, every left $\mathcal{M}$-approximation of $A$ is a monomorphism and every right $\mathcal{M}$-approximation of $A$ is an epimorphism.}
\end{remark}

A full subcategory $\mathcal{C}$ of an abelian category $\mathcal{A}$ is called {\em $n$-rigid}, if for every two objects $M,N\in \mathcal{C}$ and for every $k\in \{1,\ldots,n-1\}$, we have $\Ext_{\mathcal{A}}^k(M,N)=0$ \cite[Page 443]{Bel}. Any $n$-cluster tilting subcategory $\mathcal{M}$ of an abelian category $\mathcal{A}$ is $n$-rigid.

The following result gives a rich source of $n$-abelian categories. Note that our main result in this paper proves another direction of this theorem.
\begin{theorem}$($\cite[Theorem 3.16]{J}$)$
Let $\mathcal{A}$ be an abelian category and $\mathcal{M}$ be an $n$-cluster tilting subcategory of $\mathcal{A}$. Then $\mathcal{M}$ is an $n$-abelian category.
\end{theorem}


\section{A representation in abelian categories}

Let $\mathcal{A}$ be an abelian category. A subcategory $\mathcal{C}$ of $\mathcal{A}$ is called a {\em Serre subcategory} if for
any exact sequence
\begin{equation}
0\rightarrow A_1\rightarrow A_2\rightarrow A_3\rightarrow 0 \notag
\end{equation}
we have that $A_2\in \mathcal{C}$ if and only if $A_1\in \mathcal{C}$ and $A_3\in \mathcal{C}$. In this case we have the quotient category $\dfrac{\mathcal{A}}{\mathcal{C}}$, that is, by definition the localisation of $\mathcal{A}$ with respect to the class of all morphisms $f:X\rightarrow Y$ such that $\rm Ker(f), \rm Coker(f)\in \mathcal{C}$. $\dfrac{\mathcal{A}}{\mathcal{C}}$ is an abelian category and there exists the quotient functor
\begin{center}
$q:\mathcal{A}\rightarrow \dfrac{\mathcal{A}}{\mathcal{C}},$
\end{center}
which is exact. $q$ is the universal exact functor that annihilates $\mathcal{C}$. For more details the readers are referred to \cite{KrB,Po} or the Gabriel's thesis \cite{G}.

Let $\mathcal{C}\subseteq \mathcal{A}$ be a Serre subcategory. If the quotient functor $q:\mathcal{A}\rightarrow \dfrac{\mathcal{A}}{\mathcal{C}}$ admits a right adjoint $r:\dfrac{\mathcal{A}}{\mathcal{C}}\rightarrow \mathcal{A}$, then $\mathcal{C}$ is called a {\em localising subcategory}. The right adjoint $r$ is called the {\em section functor}, because there exists a natural isomorphism $\eta:qr\rightarrow I_{\frac{\mathcal{A}}{\mathcal{C}}}$, where $I_{\frac{\mathcal{A}}{\mathcal{C}}}:\dfrac{\mathcal{A}}{\mathcal{C}}\rightarrow \dfrac{\mathcal{A}}{\mathcal{C}}$ is the identity functor. For a Grothendieck category $\mathcal{A}$, a Serre subcategory $\mathcal{C}$ is localising if and only if $\mathcal{C}$ is closed under coproducts \cite[Proposition 6.4 of chapter 4]{Po}.

Let $\mathcal{C}$ be a full subcategory of an abelian category $\mathcal{A}$. The full subcategory
\begin{center}
$\mathcal{C}^{\bot}=\{A\in \mathcal{A} \mid \Hom_\mathcal{A}(\mathcal{C},A)=0=\Ext^1_\mathcal{A}(\mathcal{C},A)\}$
\end{center}
of $\mathcal{A}$ is called the {\em right perpendicular to $\mathcal{C}$} \cite{GL}.

\begin{proposition}$($\cite[Lemma 2.1 and Proposition 2.2]{GL}$)$\label{T3}
Let $\mathcal{C}$ be a Serre subcategory of an abelian category $\mathcal{A}$. Then the following holds:
\begin{itemize}
\item[$(i)$]
For $A\in \mathcal{A}$ and $B\in \mathcal{C}^{\bot}$, the natural homomorphism $q_{A,B}:\Hom_{\mathcal{A}}(A,B)\rightarrow \Hom_{\frac{\mathcal{A}}{\mathcal{C}}}(q(A),q(B))$ is an isomorphism.
\item[$(ii)$]
If $\mathcal{C}$ is a localising subcategory, the restriction $q:\mathcal{C}^{\bot}\rightarrow \dfrac{\mathcal{A}}{\mathcal{C}}$ is an equivalence of categories.
\end{itemize}
\end{proposition}

Let $\mathcal{M}$ be a small $n$-abelian category. We apply the above general result to $\rm{Mod}\text{-}\mathcal{M}$. Recall that $\rm{Mod}\text{-}\mathcal{M}$ is the category of all additive contravariant functors from $\mathcal{M}$ to the category of all abelian groups, it is an abelian category and that all limits and colimits are defined point-wise. Also by Yoneda's lemma, representable functors are projective and the direct sum of all representable functors $\bigoplus_{X\in \mathcal{M}}\Hom_\mathcal{M}(-,X)$, is a generator for $\rm{Mod}\text{-}\mathcal{M}$. Thus $\rm{Mod}\text{-}\mathcal{M}$ is a Grothendieck category \cite[Proposition 5.21]{Fr}. A functor $F\in \rm{Mod}\text{-}\mathcal{M}$ is called {\em finitely presented} or {\em coherent}, if there exists an exact sequence of the form
\begin{center}
$\begin{CD}
\Hom_\mathcal{M}(-,X)@ > \Hom_\mathcal{M}(-,f)>> \Hom_\mathcal{M}(-,Y)\rightarrow F\rightarrow 0.
\end{CD}$
\end{center}
We denote by $\rm{mod}\text{-}\mathcal{M}$ the full subcategory of $\rm{Mod}\text{-}\mathcal{M}$ consist of all finitely presented functors. Since every morphism in $\mathcal{M}$ has a weak kernel, $\rm{mod}\text{-}\mathcal{M}$ is an abelian category \cite[Theorem 1.4]{Fr2}.

A functor $F\in \rm{mod}\text{-}\mathcal{M}$ is called {\em effaceable}, if there is an exact sequence
\begin{center}
$\begin{CD}
\Hom_\mathcal{M}(-,X)@ > \Hom_\mathcal{M}(-,f)>> \Hom_\mathcal{M}(-,Y)\rightarrow F\rightarrow 0,
\end{CD}$
\end{center}
for some epimorphism $f:X\rightarrow Y$. The full subcategory of effaceable functors is denoted by $\rm{mod}_0\text{-}\mathcal{M}$.

\begin{definition}
Let $\mathcal{M}$ be an $n$-abelian category, $\mathcal{A}$ an abelian category and $F:\mathcal{M}\rightarrow \mathcal{A}$ a covariant additive functor.
\begin{itemize}
\item[(i)] $F$ is called {\em left $n$-exact} if for any left $n$-exact sequence $X^0 \overset{f^0}{\rightarrow} X^1 \overset{f^1}{\rightarrow} \cdots \overset{f^{n-1}}{\rightarrow} X^n \overset{f^n}{\rightarrow} X^{n+1} \notag$ in $\mathcal{M}$, $0 \rightarrow F(X^0) \rightarrow F(X^1) \rightarrow \cdots \rightarrow F(X^n)
\rightarrow F(X^{n+1}) \notag$ is an exact sequence of $\mathcal{A}$.
\item[(ii)] $F$ is called {\em right $n$-exact} if for any right $n$-exact sequence $X^0 \overset{f^0}{\rightarrow} X^1 \overset{f^1}{\rightarrow} \cdots \overset{f^{n-1}}{\rightarrow} X^n \overset{f^n}{\rightarrow} X^{n+1} \notag$ in $\mathcal{M}$, $F(X^0) \rightarrow F(X^1) \rightarrow \cdots \rightarrow F(X^n)
\rightarrow F(X^{n+1})\rightarrow 0 \notag$ is an exact sequence of $\mathcal{A}$.
\item[(iii)] $F$ is called {\em $n$-exact} if for any $n$-exact sequence $X^0 \overset{f^0}{\rightarrow} X^1 \overset{f^1}{\rightarrow} \cdots \overset{f^{n-1}}{\rightarrow} X^n \overset{f^n}{\rightarrow} X^{n+1} \notag$ in $\mathcal{M}$, $0 \rightarrow F(X^0) \rightarrow F(X^1) \rightarrow \cdots \rightarrow F(X^n)
\rightarrow F(X^{n+1})\rightarrow 0 \notag$ is an exact sequence of $\mathcal{A}$.
\end{itemize}
The contravariant left $n$-exact (resp., right $n$-exact, $n$-exact) functors are defined similarly.
\end{definition}

The above definition is a slight generalisation of \cite[Section 4.1]{Lu}. Also Bennett-Tennenhaus and Shah in \cite[Definition 2.18]{BS} defined the term "$n$-exact functor" which is different to the above definition. According to \cite[Definition 2.18]{BS} a functor $F$ between $n$-exact categories is called $n$-exact if it sends admissible $n$-exact sequences to admissible $n$-exact sequences.

\begin{proposition}$($\cite[Proposition 3.2]{EN}$)$\label{PnE}
Let $\mathcal{M}$ be an $n$-abelian category and $\mathcal{A}$ an abelian category. A covariant functor $F:\mathcal{M}\rightarrow \mathcal{A}$ is an $n$-exact functor if and only if it is both left and right $n$-exact.
\end{proposition}

In the following for simplicity we denote the representable functor $\Hom_\mathcal{M}(-,X)$ by $H_X$. The following lemma is a higher-dimensional version of Proposition 3.2 of \cite{Au}.

\begin{lemma}\label{LE}
Let $\mathcal{M}$ be a small $n$-abelian category and
$F\in \rm{mod}\text{-}\mathcal{M}$. Then the following statements are equivalent:
\begin{itemize}
\item[(1)]
$F\in \rm{mod}_0\text{-}\mathcal{M}$.
\item[(2)]
For every left $n$-exact functor $G\in \rm{mod}\text{-}\mathcal{M}$, $\Ext^{0,...,n}_{\rm{mod}\text{-}\mathcal{M}}(F,G)=0$.
\item[(3)]
For every left $n$-exact functor $G\in \rm{mod}\text{-}\mathcal{M}$, $\Hom_{\rm{mod}\text{-}\mathcal{M}}(F,G)=0$.
\item[(4)]
For every $X\in \mathcal{M}$, $\Hom_{\rm{mod}\text{-}\mathcal{M}}(F,H_X)=0$.
\end{itemize}
\begin{proof}
By Yoneda's lemma, the implication $(1)\Rightarrow (2)$ is easy. The implications $(2)\Rightarrow (3)$ and $(3)\Rightarrow (4)$ are trivial. We show the implication $(4)\Rightarrow (1)$. Consider the projective presentation
\begin{center}
$H_{Y^n}\rightarrow H_{Y^{n+1}}\rightarrow F\rightarrow 0,$
\end{center}
of $F$. Because every morphism in $\mathcal{M}$ has $n$-kernel, we have the following exact sequence.
\begin{center}
$0\rightarrow H_{Y^0}\rightarrow H_{Y^1}\rightarrow \cdots \rightarrow H_{Y^n}\rightarrow H_{Y^{n+1}}\rightarrow F\rightarrow 0.$
\end{center}
Let $X\in \mathcal{M}$. Applying $\Hom_{\rm{mod}\text{-}\mathcal{M}}(-,H_X)$ to the above exact sequence, by Yoneda's lemma, we have the exact sequence
\begin{center}
$0=\Hom_{\rm{mod}\text{-}\mathcal{M}}(F,H_X)\rightarrow \Hom_{\rm{mod}\text{-}\mathcal{M}}(H_{Y^{n+1}},H_X)\rightarrow \Hom_{\rm{mod}\text{-}\mathcal{M}}(H_{Y^n},H_X)$
\end{center}
which shows
that $Y^n\rightarrow Y^{n+1}$ is an epimorphism.
\end{proof}
\end{lemma}

\begin{proposition}
Let $\mathcal{M}$ be a small $n$-abelian category.
Then $\rm{mod}_0\text{-}\mathcal{M}$ is a Serre subcategory of $\rm{mod}\text{-}\mathcal{M}$.
\begin{proof}
It follows from Lemma \ref{LE}.
\end{proof}
\end{proposition}

We denote by $\widetilde{H}$, the composition of the Yoneda functor by the localisation functor.
\begin{center}
$\widetilde{H}:\mathcal{M}\overset{H}{\longrightarrow} \rm{mod}\text{-}\mathcal{M}\overset{q}{\longrightarrow} \dfrac{\rm{mod}\text{-}\mathcal{M}}{\rm{mod}_0\text{-}\mathcal{M}}.$
\end{center}
For simplicity we set $\mathcal{B}:=\dfrac{\rm{mod}\text{-}\mathcal{M}}{\rm{mod}_0\text{-}\mathcal{M}}$.
By Lemma \ref{LE}, for every $X\in \mathcal{M}$, $H_X\in (\rm{mod}_0\text{-}\mathcal{M})^{\bot}$ and hence by Proposition \ref{T3}(i), $\widetilde{H}$ is full and faithful. Also it is clear that $\widetilde{H}$ is an $n$-exact functor.

\section{The main theorem}
In the previous section, for a small $n$-abelian category $\mathcal{M}$, we constructed an abelian category $\mathcal{B}=\dfrac{\rm{mod}\text{-}\mathcal{M}}{\rm{mod}_0\text{-}\mathcal{M}}$ and a fully faithful functor
\begin{center}
$\widetilde{H}:\mathcal{M}\longrightarrow \mathcal{B},$
\end{center}
with the following properties:
\begin{itemize}
\item[(i)] The functor $\widetilde{H}$ is $n$-exact.
\item[(ii)] Every object $B\in \mathcal{B}$ is isomorphic to $Coker(\widetilde{H}(f))$ for some morphism $f$ in $\mathcal{M}$.
\end{itemize}
Let $\mathcal{N}$ be the essential image of $\widetilde{H}$. Obviously the induced functor $\widetilde{H}:\mathcal{M}\rightarrow \mathcal{N}$ is an equivalence. In this section we will show that $\mathcal{N}$ is an $n$-cluster tilting subcategory of the abelian category $\mathcal{B}$, which is our main result in this paper.

\begin{proposition}\label{M1}
$\mathcal{N}$ is a generating, cogenerating and functorially finite subcategory of $\mathcal{B}$.
\begin{proof}
Let $B\in \mathcal{B}$. There exists a morphism $f^0:X^0\rightarrow X^1$ in $\mathcal{M}$ such that $B$ is isomorphic to $Coker(\widetilde{H}(f^0))$. Let $X^1\overset{f^1}{\rightarrow} X^2\overset{f^2}{\rightarrow}\cdots\overset{f^{n-1}}{\rightarrow} X^n\overset{f^n}{\rightarrow} X^{n+1}$ be an $n$-cokernel of $f^0$ in $\mathcal{M}$. Then by Proposition \ref{PnE} we have the following commutative diagram with an exact row in $\mathcal{B}$.
\begin{center}
\begin{tikzpicture}
\node (X0) at (-6,1) {$\widetilde{H}(X^0)$};
\node (X1) at (-3,1) {$\widetilde{H}(X^1)$};
\node (X2) at (0,1) {$\widetilde{H}(X^2)$};
\node (X3) at (3,1) {$\cdots$};
\node (X4) at (6,1) {$\widetilde{H}(X^{n+1})$};
\node (X5) at (8,1) {$0.$};
\node (X6) at (-1.5,-0.5) {$B$};
\draw [->,thick] (X0) -- (X1) node [midway,above] {$\widetilde{H}(f^0$)};
\draw [->,thick] (X1) -- (X2) node [midway,above] {$\widetilde{H}(f^1)$};
\draw [->,thick] (X2) -- (X3) node [midway,above] {$\widetilde{H}(f^2)$};
\draw [->,thick] (X3) -- (X4) node [midway,above] {$\widetilde{H}(f^n)$};
\draw [->,thick] (X4) -- (X5) node [midway,left] {};
\draw [->>,thick] (X1) -- (X6) node [midway,below left] {$g^1$};
\draw [>->,thick] (X6) -- (X2) node [midway,below right] {$h^1$};
\end{tikzpicture}
\end{center}

The diagram immediately shows that $\mathcal{N}$ is generating-cogenerating using the morphisms $g^1$ and $h^1$. Because $\widetilde{H}(f^1)$ is a weak cokernel of $\widetilde{H}(f^0)$ in $\mathcal{N}$, and $g^1$ is an epimorphism, it is easy to see that $h^1$ is a left $\mathcal{N}$-approximation of $B$. Since $h^1:B\rightarrow \widetilde{H}(X^2)$ is the kernel of $\widetilde{H}(f^2)$, by considering $n$-kernel of $f^2$ and using a dual argument we see that $B$ has an epimorphism right $\mathcal{N}$-approximation. Therefore $\mathcal{N}$ is a generating-cogenerating functorially finite subcategory of $\mathcal{B}$.
\end{proof}
\end{proposition}

\begin{remark}\label{REx}
By the proof of Proposition \ref{M1}, every object $B\in \mathcal{B}$ fits into exact sequences
\begin{align}
&0\rightarrow B\rightarrow \widetilde{H}(X^2)\rightarrow \cdots\rightarrow \widetilde{H}(X^n)\rightarrow \widetilde{H}(X^{n+1})\rightarrow 0 \notag \\
&0\rightarrow \widetilde{H}(Y^{n+1})\rightarrow \widetilde{H}(Y^n)\rightarrow \cdots\rightarrow \widetilde{H}(Y^2)\rightarrow B\rightarrow 0 \notag
\end{align}
with $\widetilde{H}(X^2), \ldots, \widetilde{H}(X^{n+1}), \widetilde{H}(Y^2), \ldots, \widetilde{H}(Y^{n+1})\in \mathcal{N}$, such that $B\rightarrow \widetilde{H}(X^2)$ is a monomorphism left $\mathcal{N}$-approximation of $B$ and $\widetilde{H}(Y^2)\rightarrow B$ is an epimorphism right $\mathcal{N}$-approximation of $B$.
\end{remark}

\begin{lemma}\label{Ad}
Let $\mathcal{A}$ be an abelian category and $\mathcal{C}$ be a small additive subcategory of $\mathcal{A}$ which has weak kernels. Then the inclusion functor $i:\mathcal{C}\hookrightarrow \mathcal{A}$ can be extended uniquely (up to equivalence) to a right exact functor $\mathcal{F}:\rm{mod}\text{-}\mathcal{C}\rightarrow \mathcal{A}$ making the following diagram commutative, where $\mathcal{Y}$ is the Yoneda functor.
\begin{center}
\begin{tikzpicture}
\node (X0) at (0,0) {$\mathcal{C}$};
\node (X1) at (4,0) {$\rm{mod}\text{-}\mathcal{C}$};
\node (X2) at (2,-2) {$\mathcal{A}$};
\draw [->,thick] (X0) -- (X1) node [midway,above] {$\mathcal{Y}$};
\draw [->,thick] (X0) -- (X2) node [midway,below left] {$i$};
\draw [->,thick,dashed] (X1) -- (X2) node [midway,below right] {$\mathcal{F}$};
\end{tikzpicture}
\end{center}
Furthermore if $\mathcal{C}$ is also a generating and contravariantly finite subcategory of $\mathcal{A}$, then the functor $\mathcal{F}$ has a fully faithful right adjoint.
\begin{proof}
The first statement is well known, see for example \cite[Lemma 2.1.7]{KrB}.
We recall from \cite[Lemma 2.1.7]{KrB} that for every $C\in \mathcal{C}$, $\mathcal{F}(\Hom_\mathcal{C}(-, C))=C$ and for each $F\in \rm{mod}\text{-}\mathcal{C}$ with a projective presentation
\begin{center}
$\begin{CD}
\Hom_\mathcal{C}(-, C_2)@ > \Hom_\mathcal{C}(-,h)>> \Hom_\mathcal{C}(-, C_1)\longrightarrow F\longrightarrow 0,
\end{CD}$
\end{center}
we have $\mathcal{F}(F)=\Coker(h)$.

Now assume that $\mathcal{C}$ is also a generating and contravariantly finite subcategory of $\mathcal{A}$. We construct a fully faithful functor $\mathcal{R}:\mathcal{A}\rightarrow \modd$-$\mathcal{C}$ which is a right adjoint of $\mathcal{F}$. Define a functor $\mathcal{R}:\mathcal{A}\rightarrow \modd$-$\mathcal{C}$ given by $\mathcal{R}(A)=\Hom_\mathcal{A}(-, A)|_{\mathcal{C}}$ for each $A\in \mathcal{A}$. For each $A\in \mathcal{A}$, since $\mathcal{C}$ is a contravariantly finite and generating subcategory of $\mathcal{A}$, we have the following commutative diagram with an exact row
\begin{center}
\begin{tikzpicture}
\node (X0) at (3,1) {$C_2$};
\node (X1) at (5,1) {$C_1$};
\node (X2) at (7,1) {$A$};
\node (X4) at (9,1) {$0,$};
\node (X3) at (4,0) {$K$};
\draw [->,thick] (X0) -- (X1) node [midway,left] {};
\draw [->,thick] (X2) -- (X4) node [midway,left] {};
\draw [->,thick] (X1) -- (X2) node [midway,above] {f};
\draw [->,thick] (X0) -- (X3) node [midway,left] {g};
\draw [->,thick] (X3) -- (X1) node [midway,right] {h};
\end{tikzpicture}
\end{center}
where $f$ and $g$ are right $\mathcal{C}$-approximations and $h$ is the kernel of $f$ in $\mathcal{A}$. Then the sequence
\begin{center}
$\Hom_\mathcal{C}(-, C_2)\longrightarrow\Hom_\mathcal{C}(-, C_1)\longrightarrow \Hom_\mathcal{A}(-, A)|_{\mathcal{C}}\longrightarrow 0$
\end{center}
is exact and hence $\Hom_\mathcal{A}(-, A)|_{\mathcal{C}}\in \modd$-$\mathcal{C}$. Now by this exact sequence and Yoneda's lemma we have $\Hom_\mathcal{A}(A, A')\cong \Hom_{\rm{mod}\text{-}\mathcal{C}}(\Hom_\mathcal{A}(-, A)|_{\mathcal{C}}, \Hom_\mathcal{A}(-, A')|_{\mathcal{C}})$, for every $A, A'\in \mathcal{A}$. Therefore $\mathcal{R}$ is full and faithful. Now we show that $\mathcal{R}$ is a right adjoint of $\mathcal{F}$. Let $F\in \modd$-$\mathcal{C}$ and
\begin{center}
$\Hom_\mathcal{C}(-, C_2)\longrightarrow\Hom_\mathcal{C}(-, C_1)\longrightarrow F\longrightarrow 0$
\end{center}
be a projective presentation of $F$, where $C_1, C_2\in \mathcal{C}$. Since $\mathcal{F}$ is right exact, after applying $\mathcal{F}$ we have an exact sequence
$C_2\longrightarrow C_1\longrightarrow \mathcal{F}(F)\longrightarrow 0$.
Now applying $\Hom_{\rm{mod}\text{-}\mathcal{C}}(-, \Hom_\mathcal{A}(-, A)|_{\mathcal{C}})$ to the projective presentation of $F$, by Yoneda's lemma, we have an exact sequence
\begin{center}
$0\longrightarrow\Hom_{\rm{mod}\text{-}\mathcal{C}}(F, \Hom_\mathcal{A}(-, A)|_{\mathcal{C}})\longrightarrow\Hom_\mathcal{A}(C_1, A)\longrightarrow \Hom_\mathcal{A}(C_2, A)$.
\end{center}
Therefore we have $\Hom_{\rm{mod}\text{-}\mathcal{C}}(F, \Hom_\mathcal{A}(-, A)|_{\mathcal{C}})\cong \Hom_{\mathcal{A}}(\mathcal{F}(F), A)$ and the result follows.
\end{proof}
\end{lemma}

\begin{proposition}\label{P2}
Let $\mathcal{M}$ be a small $n$-abelian category. Then $\rm{mod}_0\text{-}\mathcal{M}$ is a localising subcategory of $\rm{mod}\text{-}\mathcal{M}$.
\begin{proof}
Let $\mathcal{N}$ be the essential image of $\widetilde{H}:\mathcal{M}\rightarrow \mathcal{B}=\dfrac{\rm{mod}\text{-}\mathcal{M}}{\rm{mod}_0\text{-}\mathcal{M}}$ and $\mathcal{F}:\rm{mod}\text{-}\mathcal{N}\rightarrow \mathcal{B}$ be the functor defined in Lemma \ref{Ad}.
The canonical equivalence $\mathcal{M}\rightarrow \mathcal{N}$ extends uniquely to an equivalence $\rm{mod}\text{-}\mathcal{M}\rightarrow \rm{mod}\text{-}\mathcal{N}$ making the left-hand square of the following diagram commutative.
\begin{center}
\begin{tikzpicture}
\node (X1) at (-3,0) {$\mathcal{M}$};
\node (X2) at (0,0) {$\rm{mod}\text{-}\mathcal{M}$};
\node (X3) at (3,0) {$\mathcal{B}$};
\node (X4) at (-3,-2) {$\mathcal{N}$};
\node (X5) at (0,-2) {$\rm{mod}\text{-}\mathcal{N}$};
\node (X6) at (3,-2) {$\mathcal{B}$};
\draw [->,thick] (X1) -- (X2) node [midway,left] {};
\draw [->,thick] (X2) -- (X3) node [midway,above] {$q$};
\draw [->,thick] (X1) -- (X4) node [midway,left] {$\simeq$};
\draw [->,thick] (X2) -- (X5) node [midway,left] {$\simeq$};
\draw [double,-,thick] (X3) -- (X6) node [midway,above] {};
\draw [->,thick] (X4) -- (X5) node [midway,left] {};
\draw [->,thick] (X5) -- (X6) node [midway,below] {$\mathcal{F}$};
\end{tikzpicture}
\end{center}
The commutativity of the right-hand square is an easy implication of the definition of $\mathcal{F}$. Now by Lemma \ref{Ad}, $\mathcal{F}$ and hence $q$ has a fully faithful right adjoint.
\end{proof}
\end{proposition}

The next aim is to show that $\mathcal{N}$ is an $n$-rigid subcategory of $\mathcal{B}$.\\

Let $\mathcal{A}$ be an abelian category, $k$ a positive integer and $A, C\in \mathcal{A}$. An exact sequence
\begin{equation}
0\rightarrow A\rightarrow X_{k-1}\rightarrow \cdots \rightarrow X_0\rightarrow C\rightarrow 0 \notag
\end{equation}
in $\mathcal{A}$ is called a {\em $k$-fold extension of $C$ by $A$}. Two $k$-fold extensions of $C$ by $A$,
\begin{equation}
\xi:0\rightarrow A\rightarrow B_{k-1}\rightarrow \cdots \rightarrow B_0\rightarrow C\rightarrow 0 \notag
\end{equation}
and
\begin{equation}
\xi':0\rightarrow A\rightarrow B'_{k-1}\rightarrow \cdots \rightarrow B'_0\rightarrow C\rightarrow 0 \notag
\end{equation}
are said to be {\em Yoneda equivalent} if there is a chain of $k$-fold extensions of $C$ by $A$
\begin{equation}
\xi=\xi_0,  \xi_1, \ldots ,\xi_{l-1}, \xi_l=\xi' \notag
\end{equation}
such that for every $i\in \{0,\cdots,l-1\}$, we have either a chain map $\xi_i\rightarrow \xi_{i+1}$ starting with $1_A$ and ending with $1_C$, or a chain map $\xi_{i+1}\rightarrow \xi_{i}$ starting with $1_A$ and ending with $1_C$.
$\Ext_{\mathcal{A}}^k(C,A)$ is defined as the set of Yoneda equivalence classes of $k$-fold extensions of $C$ by $A$ \cite{Mac, Mi}.

\begin{remark}\label{R}

Let $ 0\rightarrow K\rightarrow B_0\rightarrow C\rightarrow 0$ and $0\rightarrow A\rightarrow B_1'\rightarrow K'\rightarrow 0$ be two short exact sequences in an abelian category $\mathcal{A}$ and $\sigma:K\rightarrow K'$ be a morphism. Then by taking the pullback and the pushout along $\sigma$, we have the following commutative diagram with exact rows.
\begin{center}
\begin{tikzpicture}
\node (X1) at (-4.,1) {$0$};
\node (X2) at (-2.5,1) {$A$};
\node (X3) at (-1,1) {$W_1$};
\node (X4) at (0,0) {$K$};
\node (X5) at (1,1) {$B_0$};
\node (X6) at (2.5,1) {$C$};
\node (X7) at (4,1) {$0$};
\node (X8) at (-4,-1) {$0$};
\node (X9) at (-2.5,-1) {$A$};
\node (X10) at (-1,-1) {$B'_1$};
\node (X11) at (0,-2) {$K'$};
\node (X12) at (1,-1) {$W_0$};
\node (X13) at (2.5,-1) {$C$};
\node (X14) at (4,-1) {$0$};
\draw [->,thick] (X1) -- (X2) node [midway,left] {};
\draw [->,thick] (X2) -- (X3) node [midway,left] {};
\draw [->,thick] (X3) -- (X4) node [midway,left] {};
\draw [->,thick] (X3) -- (X5) node [midway,left] {};
\draw [->,thick] (X4) -- (X5) node [midway,left] {};
\draw [->,thick] (X5) -- (X6) node [midway,left] {};
\draw [->,thick] (X6) -- (X7) node [midway,above] {};
\draw [->,thick] (X8) -- (X9) node [midway,above] {};
\draw [->,thick] (X9) -- (X10) node [midway,above] {};
\draw [->,thick] (X10) -- (X11) node [midway,above] {};
\draw [->,thick] (X10) -- (X12) node [midway,left] {};
\draw [->,thick] (X11) -- (X12) node [midway,above] {};
\draw [->,thick] (X12) -- (X13) node [midway,above] {};
\draw [->,thick] (X13) -- (X14) node [midway,above] {};
\draw [double,-,thick] (X2) -- (X9) node [midway,above] {};
\draw [->,thick] (X3) -- (X10) node [midway,above] {};
\draw [->,thick] (X4) -- (X11) node [midway,above left] {$\sigma$};
\draw [->,thick] (X5) -- (X12) node [midway,above] {};
\draw [double,-,thick] (X6) -- (X13) node [midway,above] {};
\end{tikzpicture}
\end{center}
Therefore these two 2-fold extensions of $C$ by $A$ are Yoneda equivalent. Two $k$-fold extensions $S$ and $S'$ in $\mathcal{A}$ with the same start and the same end are Yoneda equivalent if and only if $S'$ can be obtained from $S$ by a finite sequence of replacements of this form and replacements of a short exact sequence by another short exact sequence equivalent to it (see \cite[Proposition 3.1 of chapter VII]{Mi}). Now let $k\geq 3$ be a positive integer and consider the following commutative diagram with exact rows
\begin{center}
\begin{tikzpicture}
\node (X4) at (0,0.5) {$0$};
\node (X5) at (1.5,0.5) {$X^0$};
\node (X6) at (3,0.5) {$X^1$};
\node (X7) at (4.5,0.5) {$\cdots$};
\node (X8) at (6,0.5) {$X^{k-2}$};
\node (X9) at (7.5,0.5) {$X^{k-1}$};
\node (X10) at (9,0.5) {$X^k$};
\node (X1) at (10.5,0.5) {$0$};
\node (X11) at (0,-1) {$\xi:0$};
\node (X12) at (1.5,-1) {$Y^0$};
\node (X13) at (3,-1) {$Y^1$};
\node (X14) at (4.5,-1) {$\cdots$};
\node (X15) at (6,-1) {$Y^{k-2}$};
\node (X16) at (7.5,-1) {$Y^{k-1}$};
\node (X17) at (9,-1) {$X^k$};
\node (X18) at (10.5,-1) {$0.$};
\draw [->,thick] (X4) -- (X5) node [midway,left] {};
\draw [->,thick] (X5) -- (X6) node [midway,left] {};
\draw [->,thick] (X10) -- (X1) node [midway,left] {};
\draw [->,thick] (X9) -- (X10) node [midway,above] {$f^{k-1}$};
\draw [->,thick] (X6) -- (X7) node [midway,above] {};
\draw [->,thick] (X7) -- (X8) node [midway,above] {};
\draw [->,thick] (X8) -- (X9) node [midway,above] {};
\draw [->,thick] (X17) -- (X18) node [midway,above] {};
\draw [->,thick] (X16) -- (X17) node [midway,below] {$g^{k-1}$};
\draw [->,thick] (X11) -- (X12) node [midway,above] {};
\draw [->,thick] (X12) -- (X13) node [midway,above] {};
\draw [->,thick] (X13) -- (X14) node [midway,above] {};
\draw [->,thick] (X14) -- (X15) node [midway,above] {};
\draw [->,thick] (X15) -- (X16) node [midway,above] {};
\draw [->,thick] (X9) -- (X16) node [midway,above] {};
\draw [double,-,thick] (X10) -- (X17) node [midway,above] {};
\end{tikzpicture}
\end{center}
By taking a pullback
\begin{center}
\begin{tikzpicture}
\node (X1) at (0,0) {$\Ker (g^{k-1})$};
\node (X2) at (0,1.5) {$\Ker (f^{k-1})$};
\node (X3) at (-2,0) {$Y^{k-2}$};
\node (X4) at (-2,1.5) {$Y_{pb}^{k-2}$};
\draw [->,thick] (X2) -- (X1) node [midway,left] {};
\draw [->,thick] (X3) -- (X1) node [midway,left] {};
\draw [->,dashed] (X4) -- (X2) node [midway,right] {};
\draw [->,dashed] (X4) -- (X3) node [midway,left] {};
\end{tikzpicture}
\end{center}
we obtain the following commutative diagram with exact rows, where the bottom row is Yoneda equivalent to $\xi$.
\begin{center}
\begin{tikzpicture}
\node (X4) at (0,0.5) {$0$};
\node (X5) at (1.5,0.5) {$X^0$};
\node (X6) at (3,0.5) {$X^1$};
\node (X7) at (4.5,0.5) {$\cdots$};
\node (X8) at (6,0.5) {$X^{k-2}$};
\node (X9) at (7.5,0.5) {$X^{k-1}$};
\node (X10) at (9,0.5) {$X^k$};
\node (X1) at (10.5,0.5) {$0$};
\node (X11) at (0,-1) {$0$};
\node (X12) at (1.5,-1) {$Y^0$};
\node (X13) at (3,-1) {$Y^1$};
\node (X14) at (4.5,-1) {$\cdots$};
\node (X15) at (6,-1) {$Y_{pb}^{k-2}$};
\node (X16) at (7.5,-1) {$X^{k-1}$};
\node (X17) at (9,-1) {$X^k$};
\node (X18) at (10.5,-1) {$0.$};
\draw [->,thick] (X4) -- (X5) node [midway,left] {};
\draw [->,thick] (X5) -- (X6) node [midway,left] {};
\draw [->,thick] (X10) -- (X1) node [midway,left] {};
\draw [->,thick] (X9) -- (X10) node [midway,above] {};
\draw [->,thick] (X6) -- (X7) node [midway,above] {};
\draw [->,thick] (X7) -- (X8) node [midway,above] {};
\draw [->,thick] (X8) -- (X9) node [midway,above] {};
\draw [->,thick] (X17) -- (X18) node [midway,above] {};
\draw [->,thick] (X16) -- (X17) node [midway,below] {};
\draw [->,thick] (X11) -- (X12) node [midway,above] {};
\draw [->,thick] (X12) -- (X13) node [midway,above] {};
\draw [->,thick] (X13) -- (X14) node [midway,above] {};
\draw [->,thick] (X14) -- (X15) node [midway,above] {};
\draw [->,thick] (X15) -- (X16) node [midway,above] {};
\draw [double,-,thick] (X9) -- (X16) node [midway,above] {};
\draw [double,-,thick] (X10) -- (X17) node [midway,above] {};
\end{tikzpicture}
\end{center}
\end{remark}

Now we can prove that $\mathcal{N}$ is $n$-rigid.

\begin{theorem}\label{TR}
$\mathcal{N}$ is an $n$-rigid subcategory of the abelian category $\mathcal{B}$.
\begin{proof}
Let $k\in \{1,\ldots,n-1\}$, $N_1, N_2\in \mathcal{N}$ and
\begin{center}
$\xi: 0\rightarrow N_2\rightarrow F^{n-k+1}\rightarrow \cdots \rightarrow F^n\rightarrow N_1\rightarrow 0$
\end{center}
be a $k$-fold extension of $N_1$ by $N_2$ in $\mathcal{B}$. There exists $X',X^{n+1}\in \mathcal{M}$ such that $\widetilde{H}(X^{n+1})=N_1$ and $\widetilde{H}(X')=N_2$. The cokernel of $F^n\rightarrow H_{X^{n+1}}$ in $\rm{mod}\text{-}\mathcal{M}$, denoted by $C$, is effaceable.
Since $F^n$ is a finitely presented functor, there is an epimorphism $H_{X^n}\rightarrow F^n$ in $\rm{mod}\text{-}\mathcal{M}$ for some $X^n\in \mathcal{M}$. Obviously the cokernel of the composition $H_{X^n}\rightarrow F^n\rightarrow H_{X^{n+1}}$ is equal to $C$. By Yoneda's lemma this composition is induced by a morphism $f:X^n\rightarrow X^{n+1}$. The proof of the implication $(4)\Rightarrow (1)$ of Lemma \ref{LE} shows that $f$ is an epimorphism. Since $\mathcal{M}$ is $n$-abelian, $f$ sits in an $n$-exact sequence $0\rightarrow X^0\rightarrow X^1\rightarrow \cdots \rightarrow X^n\overset{f}{\rightarrow} X^{n+1}\rightarrow 0.$
Then we have the following commutative diagram with exact rows in $\mathcal{B}$.
\begin{center}
\begin{tikzpicture}
\node (X1) at (-10,1) {$0$};
\node (X2) at (-8.5,1) {$\widetilde{H}(X^0)$};
\node (X3) at (-7,1) {$\cdots$};
\node (X4) at (-5,1) {$\widetilde{H}(X^{n-k})$};
\node (X5) at (-2,1) {$\widetilde{H}(X^{n-k+1})$};
\node (X6) at (0,1) {$\cdots$};
\node (X7) at (2,1) {$\widetilde{H}(X^n)$};
\node (X8) at (4,1) {$\widetilde{H}(X^{n+1})$};
\node (X9) at (5.75,1) {$0$};
\node (X10) at (-6,-0.5) {$0$};
\node (X11) at (-4,-0.5) {$\widetilde{H}(X')$};
\node (X12) at (-2,-0.5) {$F^{n-k+1}$};
\node (X13) at (0,-0.5) {$\cdots$};
\node (X14) at (2,-0.5) {$F^n$};
\node (X15) at (4,-0.5) {$\widetilde{H}(X^{n+1})$};
\node (X16) at (5.75,-0.5) {$0.$};
\draw [->,thick] (X1) -- (X2) node [midway,left] {};
\draw [->,thick] (X2) -- (X3) node [midway,left] {};
\draw [->,thick] (X3) -- (X4) node [midway,left] {};
\draw [->,thick] (X4) -- (X5) node [midway,left] {};
\draw [->,thick] (X5) -- (X6) node [midway,left] {};
\draw [->,thick] (X6) -- (X7) node [midway,above] {};
\draw [->,thick] (X7) -- (X8) node [midway,left] {};
\draw [->,thick] (X8) -- (X9) node [midway,above] {};
\draw [->,thick] (X10) -- (X11) node [midway,left] {};
\draw [->,thick] (X11) -- (X12) node [midway,above] {};
\draw [->,thick] (X12) -- (X13) node [midway,above] {};
\draw [->,thick] (X13) -- (X14) node [midway,above] {};
\draw [->,thick] (X14) -- (X15) node [midway,above] {};
\draw [->,thick] (X15) -- (X16) node [midway,above] {};
\draw [->,thick] (X7) -- (X14) node [midway,above] {};
\draw [double,-,thick] (X8) -- (X15) node [midway,above] {};
\end{tikzpicture}
\end{center}
First assume that $k=1$. By the universal property of kernel there exists a morphism $\widetilde{H}(X^{n-1})\rightarrow \widetilde{H}(X')$ that makes the solid part of the following diagram commutative.
\begin{center}
\begin{tikzpicture}
\node (X3) at (-9,0.5) {$0$};
\node (X4) at (-7,0.5) {$\widetilde{H}(X^0$)};
\node (X5) at (-5.5,0.5) {$\cdots$};
\node (X6) at (-3.5,0.5) {$\widetilde{H}(X^{n-1})$};
\node (X7) at (-1.5,0.5) {$\widetilde{H}(X^n)$};
\node (X8) at (0.5,0.5) {$\widetilde{H}(X^{n+1})$};
\node (X9) at (2,0.5) {$0$};
\node (X10) at (-5,-1) {$0$};
\node (X11) at (-3.5,-1) {$\widetilde{H}(X')$};
\node (X12) at (-1.5,-1) {$F^n$};
\node (X13) at (0.5,-1) {$\widetilde{H}(X^{n+1})$};
\node (X14) at (2,-1) {$0.$};
\draw [->,thick] (X3) -- (X4) node [midway,left] {};
\draw [->,thick] (X4) -- (X5) node [midway,left] {};
\draw [->,thick] (X5) -- (X6) node [midway,left] {};
\draw [->,thick] (X6) -- (X7) node [midway,above] {};
\draw [->,thick] (X7) -- (X8) node [midway,left] {};
\draw [->,thick] (X8) -- (X9) node [midway,above] {};
\draw [->,thick] (X10) -- (X11) node [midway,left] {};
\draw [->,thick] (X11) -- (X12) node [midway,above] {};
\draw [->,thick] (X12) -- (X13) node [midway,above] {};
\draw [->,thick] (X13) -- (X14) node [midway,above] {};
\draw [->,thick] (X6) -- (X11) node [midway,above] {};
\draw [->,thick] (X7) -- (X12) node [midway,above] {};
\draw [double,-,thick] (X8) -- (X13) node [midway,above] {};
\draw [->,thick,dashed] (X7) -- (X11) node [midway,above] {};
\end{tikzpicture}
\end{center}
Since the top row is induced by an $n$-exact sequence, in particular $\widetilde{H}(X^{n-1})\rightarrow \widetilde{H}(X^n)$ is a weak cokernel of $\widetilde{H}(X^{n-2})\rightarrow \widetilde{H}(X^{n-1})$ in $\mathcal{N}$, and the morphism $\widetilde{H}(X')\rightarrow F^n$ is a monomorphism, there exists a morphism $\widetilde{H}(X^n)\dashrightarrow \widetilde{H}(X')$ that makes the triangle containing $\widetilde{H}(X^{n-1})$, $\widetilde{H}(X^n)$ and $\widetilde{H}(X')$ commutative. By \cite[Chapter 1, Proposition 5.6]{ARS} the left-hand square is a pushout square. The universal property of pushout implies that $\widetilde{H}(X')\rightarrow F^n$ has a left inverse. Then the bottom row splits. This shows that $\Ext^1_\mathcal{B}(\mathcal{N}, \mathcal{N})=0$.

Now assume that $k>1$. Using the argument in Remark \ref{R}, we have the commutative diagram
\begin{diagram}\label{d1}
\begin{center}
\begin{tikzpicture}
\node (X1) at (-10,1) {$0$};
\node (X2) at (-8.5,1) {$\widetilde{H}(X^0)$};
\node (X3) at (-7,1) {$\cdots$};
\node (X4) at (-5,1) {$\widetilde{H}(X^{n-k})$};
\node (X5) at (-2,1) {$\cdots$};
\node (X6) at (0,1) {$\widetilde{H}(X^{n-1})$};
\node (X7) at (2,1) {$\widetilde{H}(X^n)$};
\node (X8) at (4,1) {$\widetilde{H}(X^{n+1})$};
\node (X9) at (5.75,1) {$0$};
\node (X10) at (-6,-0.5) {$0$};
\node (X11) at (-4,-0.5) {$\widetilde{H}(X')$};
\node (X12) at (-2,-0.5) {$\cdots$};
\node (X13) at (0,-0.5) {$F^{n-1}_{pb}$};
\node (X14) at (2,-0.5) {$\widetilde{H}(X^{n})$};
\node (X15) at (4,-0.5) {$\widetilde{H}(X^{n+1})$};
\node (X16) at (5.75,-0.5) {$0,$};
\draw [->,thick] (X1) -- (X2) node [midway,left] {};
\draw [->,thick] (X2) -- (X3) node [midway,left] {};
\draw [->,thick] (X3) -- (X4) node [midway,left] {};
\draw [->,thick] (X4) -- (X5) node [midway,left] {};
\draw [->,thick] (X5) -- (X6) node [midway,left] {};
\draw [->,thick] (X6) -- (X7) node [midway,above] {};
\draw [->,thick] (X7) -- (X8) node [midway,left] {};
\draw [->,thick] (X8) -- (X9) node [midway,above] {};
\draw [->,thick] (X10) -- (X11) node [midway,left] {};
\draw [->,thick] (X11) -- (X12) node [midway,above] {};
\draw [->,thick] (X12) -- (X13) node [midway,above] {};
\draw [->,thick] (X13) -- (X14) node [midway,above] {};
\draw [->,thick] (X14) -- (X15) node [midway,above] {};
\draw [->,thick] (X15) -- (X16) node [midway,above] {};
\draw [double,-,thick] (X7) -- (X14) node [midway,above] {};
\draw [double,-,thick] (X8) -- (X15) node [midway,above] {};
\end{tikzpicture}
\end{center}
\end{diagram}
with exact rows in $\mathcal{B}$, where the bottom row is Yoneda equivalent to $\xi$.
Now for $j\in\{1,\cdots,n-1\}$ we set $C^j:=\Imm(\widetilde{H}(X^j)\rightarrow \widetilde{H}(X^{j+1}))$ in $\mathcal{B}$. Indeed we split the first row of Diagram 4.1 into short exact sequences as follow.
\begin{center}
\begin{tikzpicture}
\node (X1) at (-9,0) {$0$};
\node (X2) at (-7.5,0) {$\widetilde{H}(X^0)$};
\node (X3) at (-5.5,0) {$\widetilde{H}(X^1)$};
\node (X4) at (-3.5,0) {$\widetilde{H}(X^2)$};
\node (X5) at (-1.5,0) {$\cdots$};
\node (X6) at (.5,0) {$\widetilde{H}(X^{n-1})$};
\node (X7) at (2.5,0) {$\widetilde{H}(X^n)$};
\node (X8) at (4.5,0) {$\widetilde{H}(X^{n+1})$};
\node (X9) at (6.5,0) {$0.$};
\node (X10) at (-4.5,-1) {$C^1$};
\node (X11) at (-2.5,-1) {$C^2$};
\node (X12) at (-.5,-1) {$C^{n-2}$};
\node (X13) at (1.5,-1) {$C^{n-1}$};
\node (X14) at (-1.5,-1) {$\cdots$};
\draw [->,thick] (X1) -- (X2) node [midway,left] {};
\draw [->,thick] (X2) -- (X3) node [midway,left] {};
\draw [->,thick] (X3) -- (X4) node [midway,left] {};
\draw [->,thick] (X4) -- (X5) node [midway,left] {};
\draw [->,thick] (X5) -- (X6) node [midway,left] {};
\draw [->,thick] (X6) -- (X7) node [midway,above] {};
\draw [->,thick] (X7) -- (X8) node [midway,left] {};
\draw [->,thick] (X8) -- (X9) node [midway,above] {};
\draw [->>,thick] (X3) -- (X10) node [midway,left] {};
\draw [>->,thick] (X10) -- (X4) node [midway,above] {};
\draw [->>,thick] (X4) -- (X11) node [midway,above] {};
\draw [>->,thick] (X11) -- (X5) node [midway,above] {};
\draw [->>,thick] (X5) -- (X12) node [midway,above] {};
\draw [>->,thick] (X12) -- (X6) node [midway,above] {};
\draw [->>,thick] (X6) -- (X13) node [midway,above] {};
\draw [>->,thick] (X13) -- (X7) node [midway,above] {};
\end{tikzpicture}
\end{center}
First for $j\in\{1,\cdots,n-2\}$ by applying the functor $\Hom_{\mathcal{B}}(-,\widetilde{H}(X'))$ to the short exact sequence
\begin{equation}
0\rightarrow C^j\rightarrow \widetilde{H}(X^{j+1})\rightarrow C^{j+1}\rightarrow 0, \notag
\end{equation}
we have the exact sequence

$
\Hom_{\mathcal{B}}(\widetilde{H}(X^{j+1}),\widetilde{H}(X'))\rightarrow \Hom_{\mathcal{B}}(C^j,\widetilde{H}(X'))\rightarrow\Ext^1_{\mathcal{B}}(C^{j+1},\widetilde{H}(X'))\rightarrow \\ \Ext^1_{\mathcal{B}}(\widetilde{H}(X^{j+1}),\widetilde{H}(X')). \notag
$
By the first part of the proof, $\Ext^1_\mathcal{B}(\mathcal{N}, \mathcal{N})=0$ and hence\\ $\Ext^1_{\mathcal{B}}(\widetilde{H}(X^{j+1}),\widetilde{H}(X'))=0$. Since $\widetilde{H}(X^j)\rightarrow \widetilde{H}(X^{j+1})$ is a weak cokernel of $\widetilde{H}(X^{j-1})\rightarrow \widetilde{H}(X^j)$ in $\mathcal{N}$ and $\widetilde{H}(X^j)\rightarrow C^j$ is an epimorphism, $C^j\rightarrow \widetilde{H}(X^{j+1})$ is a left $\mathcal{N}$-approximation, and so the map $\Hom_{\mathcal{B}}(\widetilde{H}(X^{j+1}),\widetilde{H}(X'))\rightarrow \Hom_{\mathcal{B}}(C^j,\widetilde{H}(X'))$ is an epimorphism. Thus $\Ext_{\mathcal{B}}^1(C^{j+1},\widetilde{H}(X'))=0$ for $j\in\{1,\cdots,n-2\}$.

Now by induction assume that $k\in\{2,\cdots,n-1\}$ and $\mathcal{N}$ is $k$-rigid. We show that $\mathcal{N}$ is $(k+1)$-rigid.
Choose an arbitrary element $\xi\in \Ext_{\mathcal{B}}^k(\widetilde{H}(X^{n+1}),\widetilde{H}(X'))$ and construct Diagram 4.1 for it. Since we have the following commutative diagram with exact row,
\begin{center}
\begin{tikzpicture}
\node (X1) at (-5,0) {$0$};
\node (X4) at (-3.5,0) {$\widetilde{H}(X')$};
\node (X5) at (-1.5,0) {$\cdots$};
\node (X6) at (.5,0) {$F^{n-1}_{pb}$};
\node (X7) at (2.5,0) {$\widetilde{H}(X^n)$};
\node (X8) at (5,0) {$\widetilde{H}(X^{n+1})$};
\node (X9) at (7,0) {$0,$};
\node (X13) at (1.5,-1) {$C^{n-1}$};
\draw [->,thick] (X1) -- (X4) node [midway,left] {};
\draw [->,thick] (X4) -- (X5) node [midway,left] {};
\draw [->,thick] (X5) -- (X6) node [midway,left] {};
\draw [->,thick] (X6) -- (X7) node [midway,above] {};
\draw [->,thick] (X7) -- (X8) node [midway,left] {};
\draw [->,thick] (X8) -- (X9) node [midway,above] {};
\draw [->>,thick] (X6) -- (X13) node [midway,above] {};
\draw [>->,thick] (X13) -- (X7) node [midway,above] {};
\end{tikzpicture}
\end{center}
it is enough to show that $\Ext_{\mathcal{B}}^{k-1}(C^{n-1},\widetilde{H}(X'))=0$ (see \cite[Page 175]{Mi}). Because by hypothesis $\Ext_{\mathcal{B}}^{1,\cdots,k-1}(\widetilde{H}(X^i),\widetilde{H}(X'))=0$ for $i\in\{0,\cdots,n+1\}$, we have
$
\Ext_{\mathcal{B}}^{k-1}(C^{n-1},\widetilde{H}(X'))\\\cong \Ext_{\mathcal{B}}^{k-2}(C^{n-2},\widetilde{H}(X'))\cong \cdots \cong\Ext_{\mathcal{B}}^{2}(C^{n-k+2},\widetilde{H}(X'))\cong \Ext_{\mathcal{B}}^{1}(C^{n-k+1},\widetilde{H}(X'))=0 \notag
$
and the result follows.
\end{proof}
\end{theorem}

Now we can prove our main result in this paper.

\begin{theorem}\label{m2}
Let $\mathcal{M}$ be a small $n$-abelian category. Then $\mathcal{N}$ is an $n$-cluster tilting subcategory of the abelian category $\mathcal{B}$.
\begin{proof}
By Proposition \ref{M1} and Theorem \ref{TR}, it is enough to show that $\mathcal{N}$ is maximal $n$-rigid. Let $B$ be an object of $\mathcal{B}$ such that for every $N\in \mathcal{N}$ and every $k\in \{1,2,\ldots,n-1\}$ we have $\Ext_{\mathcal{B}}^k(N,B)=0$. We need to show that
$B\in \mathcal{N}$. Let
\begin{equation}
0\rightarrow B\rightarrow \widetilde{H}(X^2)\rightarrow \widetilde{H}(X^3)\rightarrow \cdots \rightarrow \widetilde{H}(X^n)\rightarrow \widetilde{H}(X^{n+1})\rightarrow 0, \notag
\end{equation}
be the first exact sequence in Remark \ref{REx} with $\widetilde{H}(X^2), \cdots, \widetilde{H}(X^{n+1})\in \mathcal{N}$. Set $C^1=B$ and $C^r=\Imm(\widetilde{H}(X^r)\rightarrow \widetilde{H}(X^{r+1}))$ in $\mathcal{B}$ for $r\in \{2,\ldots n-1\}$. Inductively we show that for every $r\in \{1,2,\ldots n-1\}$, $k\in \{1,2,\ldots n-r\}$ and every $N\in \mathcal{N}$, $\Ext_{\mathcal{B}}^k(N,C^r)=0$. If $r=1$, then the claim follows from assumption. Suppose that $\Ext_{\mathcal{B}}^k(N,C^r)=0$ for each $k\in \{1,2,\ldots n-r\}$. We show that for each $k\in \{1,2,\ldots n-r-1\}$, $\Ext_{\mathcal{B}}^k(N,C^{r+1})=0$. Applying the functor $\Hom_\mathcal{B}(N,-)$ to the exact sequence $0\rightarrow C^r\rightarrow \widetilde{H}(X^{r+1})\rightarrow C^{r+1}\rightarrow 0$, gives an exact sequence
\begin{equation}
\Ext_{\mathcal{B}}^k(N,\widetilde{H}(X^{r+1}))\rightarrow \Ext_{\mathcal{B}}^k(N,C^{r+1})\rightarrow \Ext_{\mathcal{B}}^{k+1}(N,C^r). \notag
\end{equation}
By Theorem \ref{TR}, $\Ext_{\mathcal{B}}^k(N,\widetilde{H}(X^{r+1}))=0$ and hence the claim follows. Then the sequence $0\rightarrow C^{n-1}\rightarrow \widetilde{H}(X^n) \rightarrow \widetilde{H}(X^{n+1})\rightarrow 0$ is an split short exact sequence and since $\mathcal{N}$ is idempotent complete, $C^{n-1}\in \mathcal{N}$. Again because $\Ext_{\mathcal{B}}^1(C^{n-1},C^{n-2})=0$, we have $C^{n-2}\in \mathcal{N}$. Repeating this process we obtain that $C^2\in \mathcal{N}$. Since $\Ext_{\mathcal{B}}^1(C^2,B)=0$, $B$ is a direct summand of $\widetilde{H}(X^2)$ and so $B\in\mathcal{N}$.

Now let $B\in \mathcal{B}$ such that $\Ext_{\mathcal{B}}^k(B,N)=0$ for every $N\in \mathcal{N}$ and every $k\in \{1,2,\ldots,n-1\}$. Let
\begin{align}
&0\rightarrow \widetilde{H}(Y^{n+1})\rightarrow \widetilde{H}(Y^n)\rightarrow \cdots\rightarrow \widetilde{H}(Y^2)\rightarrow B\rightarrow 0, \notag
\end{align}
be the second exact sequence in Remark \ref{REx} with $\widetilde{H}(Y^2), \cdots, \widetilde{H}(Y^{n+1})\in \mathcal{N}$. Set $C^1=B$, $C^2=\Ker(\widetilde{H}(Y^2)\rightarrow B)$ and $C^r=\Ker(\widetilde{H}(Y^r)\rightarrow \widetilde{H}(Y^{r-1}))$ in $\mathcal{B}$ for $r\in \{3,\ldots n-1\}$. A similar (completely dual) argument as in the first part of the proof shows that $B \in\mathcal{N}$. Therefore $\mathcal{N}$ is an $n$-cluster tilting subcategory of $\mathcal{B}$.
\end{proof}
\end{theorem}

\section*{acknowledgements}
The authors would like to thank the referee for a careful reading of this
paper and making many helpful suggestions that improved the paper. During the preparation of this paper Sondre Kvamme informed us that he has a different proof that any $n$-abelian category is equivalent to an $n$-cluster tilting subcategory of an abelian category, but he hasn't written up the result yet (see \cite{Kv1}). The research of the second
author was in part supported by a grant from IPM (No. 99170412).

\end{document}